\documentclass[a4paper]{amsart}

\usepackage{epic, eepic, amsfonts, latexsym, amssymb, graphicx,
multicol, mathrsfs, color, amscd, verbatim, paralist,
xspace, url, euscript, stmaryrd,  amsmath, enumitem,
bbold, multirow, tikz, mathtools}
\usepackage[all,pdf,cmtip]{xy}

\usepackage[colorlinks, linkcolor=blue, citecolor=magenta, urlcolor=cyan]{hyperref}

\usepackage{graphicx}

\usepackage{float}

\def\hat{\widehat}

\def\RR{{\mathbb R}}
\def\CC{{\mathbb C}}

\def\NN{{\mathbb N}}

\newcommand{\suchthat}{\;\ifnum\currentgrouptype=16 \middle\fi|\;}

\def\ra{\rightarrow}

\def\hat{\widehat}

\def\aff{\mathop{\rm aff}\nolimits}

\def\qed{{\hfill $\Box$}}
\newtheorem{theorem}{THEOREM}[section]

\theoremstyle{definition}

\theoremstyle{remark}
\newtheorem{remark}[theorem]{Remark}

\makeatletter
\def\blfootnote{\xdef\@thefnmark{}\@footnotetext}
\makeatother

\begin{document}

\title[Hyperbolicity of tube domains in $\CC^2$]{On necessary and sufficient conditions
\vspace{0.1cm}\\
for the Kobayashi hyperbolicity
\vspace{0.1cm}\\
of tube domains in $\CC^2$}
\author[Isaev]{Alexander Isaev}

\address{Mathematical Sciences Institute\\
Australian National University\\
Canberra, Acton, ACT 2601, Australia}
\email{alexander.isaev@anu.edu.au}

\maketitle

\thispagestyle{empty}

\pagestyle{myheadings}

\begin{abstract} This note concerns tube domains in $\CC^2$ with the envelope of holomorphy not equal to the entire space. We construct examples showing that for such domains the sufficient condition for Kobayashi hyperbolicity due to M.~Jarnicki and P.~Pflug cannot be replaced by its weaker \lq\lq affine\rq\rq\, variant, which is known to be a necessary condition for hyperbolicity. Thus, we arrive at the somewhat unexpected conclusion that the obstructions for a domain in the above class to be Kobayashi hyperbolic are not just \lq\lq affine\rq\rq.  
\end{abstract}

\section{Introduction}\label{intro}
\setcounter{equation}{0}

A connected complex manifold $X$ is said to be {\it Kobayashi-hyperbolic}\, (or simply {\it hyperbolic}) if the Koba\-ya\-shi pseudodistance on $X$ is in fact a distance (see \cite{K1}, \cite{K2} for details). If $X$ is endowed with a Riemannian metric, hyperbolicity is equivalent to the following property: for any point $x\in X$ there exist a neighborhood $U$ of $x$ and a constant $M>0$ such that for all holomorphic maps $f:\Delta\ra X$ with $f(0)\in U$ one has $||df(0)||<M$, where $\Delta$ is the unit disk in $\CC$ (see, e.g., \cite{L} and \cite{HK}). Verification of hyperbolicity may be quite hard even for very special classes of manifolds.

This short note is a follow-up to our earlier paper \cite{I}. We discuss {\it tube domains}\, in $\CC^n$, i.e., domains of the form $T_D:=D+i\RR^n$, where $D$ is a domain in $\RR^n$ called the {\it base}\, of $T_D$. Clearly, for a tube domain $T_D\subset\CC^n$ hyperbolicity is equivalent to the following condition: for every point $x\in D$ there exist a neighborhood $U$ of $x$ in $D$ and a constant $M>0$ such that for all harmonic maps $f:\Delta\ra D$ with $f(0)\in U$ one has $||df(0)||<M$ (cf.~\cite{L} and \cite[Theorem 13.6.2]{JP}).

We assume that $n=2$. Surprisingly, so far no easily verifiable criterion for the hyperbolicity of a tube domain has been found even in this situation. By Bochner's theorem, the envelope of holomorphy of $T_D$ coincides with $T_{\hat D}$, where $\hat D$ is the convex hull of $D$ (see, e.g., \cite[Section 21]{V}), and it is natural to investigate hyperbolicity separately in each of the cases: (i) $T_{\hat D}\not=\CC^2$ and (ii) $T_{\hat D}=\CC^2$. In \cite{HI} we looked at several classes of hyperbolic domains in $\CC^2$ falling in case (ii). For example, we showed that $T_D$ is hyperbolic  if $D$ is a domain bounded by two spirals, where a spiral is a curve defined in polar coordinates in $\RR^2$ by the equation $r=g(\varphi)$, with $g$ being an increasing function of $\varphi$ such that $\lim_{\varphi\ra-\infty}g(\varphi)=0$ and $\lim_{\varphi\ra+\infty}g(\varphi)=\infty$. However, there is no comprehensive description of all hyperbolic domains covered by case (ii) (cf.~\cite[p.~533, Question 13.6]{JP}). 

On the other hand, for domains in $\CC^2$ falling in case (i) certain progress towards finding a hyperbolicity criterion has been made. Firstly, such a criterion was proposed by J.-J. Loeb in \cite[Th\'eor\`eme 6]{L}.  To state Loeb's result, let $D\subset\RR^2$ be a domain with $\hat D\ne\RR^2$. Writing coordinates in $\CC^2$ as $z_j=x_j+iy_j$, $j=1,2$, we may assume without loss of generality that $D$ lies in the half-space $\{x_2>0\}$. In this situation, we say that $D$ has Property (L) if the following holds:

\begin{itemize}

\item[{}] \emph{there does not exist a point $a=(a_1,a_2)\in D$ for which there is a sequence of real numbers $\{b_k\}$ converging to $a_2$ such that the segment $[-k,k]\times\{b_k\}$ lies in $D$ for all $k\in\NN$.}

\end{itemize}
\vspace{0.05cm}

\noindent Then \cite[Th\'eor\`eme 6]{L} asserts that $T_D$ is hyperbolic if and only if $D$ has Property (L). The necessity implication is obvious, but, unfortunately, the nice argument by Loeb contains a flaw, and in \cite{I} we were able to constructed counterexamples to the sufficiency implication (cf.~\cite[Part (a) of Remark 13.6.7]{JP}).

In fact, as M.~Jarnicki and P.~Pflug observed in \cite[Part (b) of Theorem 13.6.6] {JP}, the proof provided in \cite{L} only yields a weaker statement. To formulate it, we say that a domain $D$ lying in the half-space $\{x_2>0\}$ has Property (J-P) if

\begin{itemize}

\item[{}] \emph{there does not exist a point $a=(a_1,a_2)\in D$ such that for every $k\in\NN$ one can find a real-analytic function $\gamma_k(t)$ on $[-k, k]$, with $(t,\gamma_k(t))\in D$ and $|\gamma_k(t)-a_2|\le 1/k$ for all $t$.} 

\end{itemize}
\vspace{0.05cm}

\noindent Clearly, we have the implication Property (J-P)$\Rightarrow$Property (L). The result of Jarnicki and Pflug can now be stated as follows:

\begin{theorem}\label{correctedthm} Let $D\subset\RR^2$ be a domain lying in the half-space $\{x_2>0\}$ that has Property {\rm (J-P)}. Then the tube domain $T_D$ is hyperbolic.
\end{theorem}

Next, we introduce the \lq\lq affine\rq\rq\, variant of Property (J-P) (cf.~\cite{I}). Namely, we say that a domain $D$ lying in the half-space $\{x_2>0\}$ has Property (J-P$)_{\aff}$ if 

\begin{itemize}

\item[{}] \emph{there does not exist a point $a=(a_1,a_2)\in D$ such that for every $k\in\NN$ one can find an affine function $\gamma_k(t)=c_kt+d_k$, with $(t,\gamma_k(t))\in D$ and $|\gamma_k(t)-a_2|\le 1/k$ for all $t\in [-k,k]$.}

\end{itemize}

\noindent We have the implications Property (J-P)$\Rightarrow$Property (J-P$)_{\aff}$$\Rightarrow$Property (L). It is not hard to see that Property (J-P$)_{\aff}$ is a necessary condition for the hyperbolicity of $T_D$ (see \cite[Theorem 1.3]{I}).

The introduction of Property (J-P$)_{\aff}$ in \cite{I} reflected the expectation, which has been around for some time now, that the obstructions for the hyperbolicity of a tube domain $T_D$ with $T_{\hat D}\not=\CC^2$ should be \lq\lq affine\rq\rq. The examples given in \cite{I} show that Property (L) does not describe all the obstructions, so the only other natural \lq\lq affine\rq\rq\, condition appears to be the stronger Property (J-P$)_{\aff}$. In this note we demonstrate that Property (J-P$)_{\aff}$ is not sufficient for hyperbolicity either. Namely, we strengthen \cite[Theorem 1.2]{I} as follows:

\begin{theorem}\label{result} There exists a domain $D\subset\RR^2$ lying in the half-space $\{x_2>0\}$ that has Property {\rm (J-P}$)_{\aff}$ and for which $T_D$ is not hyperbolic. Such a domain $D$ can be chosen to have a $C^{\infty}$-smooth boundary.
\end{theorem}

Thus, hyperbolicity and the three properties discussed above are related as shown in the following diagram:

$$\xymatrix{
\hbox{Property (J-P)}  \ar@{=>}[r] \ar@{=>}[d] & \hbox{\rm Property (J-P$)_{\aff}$} \hspace{-0.29cm}\ar@{=>}[d] \hspace{-0.3cm}\ar@<1ex>@{=>}[ld] \\
\hspace{-0.7cm}\hbox{Hyperbolicity} \hspace{0.2cm}\ar@<0.5ex>@{=>}[r]_{\backslash} 	\ar@<0.5ex>@{=>}[ru]_{\backslash}	& \hbox{Property (L).} \ar@<1ex>@{=>}[l]
}
$$

Although the examples provided below are elementary, they are nevertheless intriguing as they suggest that the search for an \lq\lq affine\rq\rq\, criterion for hyperbolicity should be abandoned. The problem of eliminating the gap between necessary and sufficient conditions remains open, but in order to solve this problem one must look beyond \lq\lq affine\rq\rq\, properties. The next most natural question to address is then whether Property (J-P) is a necessary condition for hyperbolicity.

{\bf Acknowledgement.} This work is supported by the Australian Research\linebreak Council, grant DP140100296.

\section{The examples}\label{sect1}
\setcounter{equation}{0}

First, let
$$
\begin{array}{l}
D:=\{(x_1,x_2)\in\RR^2\mid 0<x_2<4\}\setminus\\
\vspace{-0.05cm}\\
\hspace{1cm}\Bigl(\{-3\pi/2\}\times[0,2]\cup\{-\pi/2\}\times[2,4]\cup\{\pi/2\}\times[0,2]\cup\{3\pi/2\}\times[2,4]\Bigr),
\end{array}
$$
as shown in Fig.~1 below. Clearly, $D$ has Property (J-P$)_{\aff}$.

\begin{figure}[H]
\begin{center}  
\includegraphics[width=4.2in]{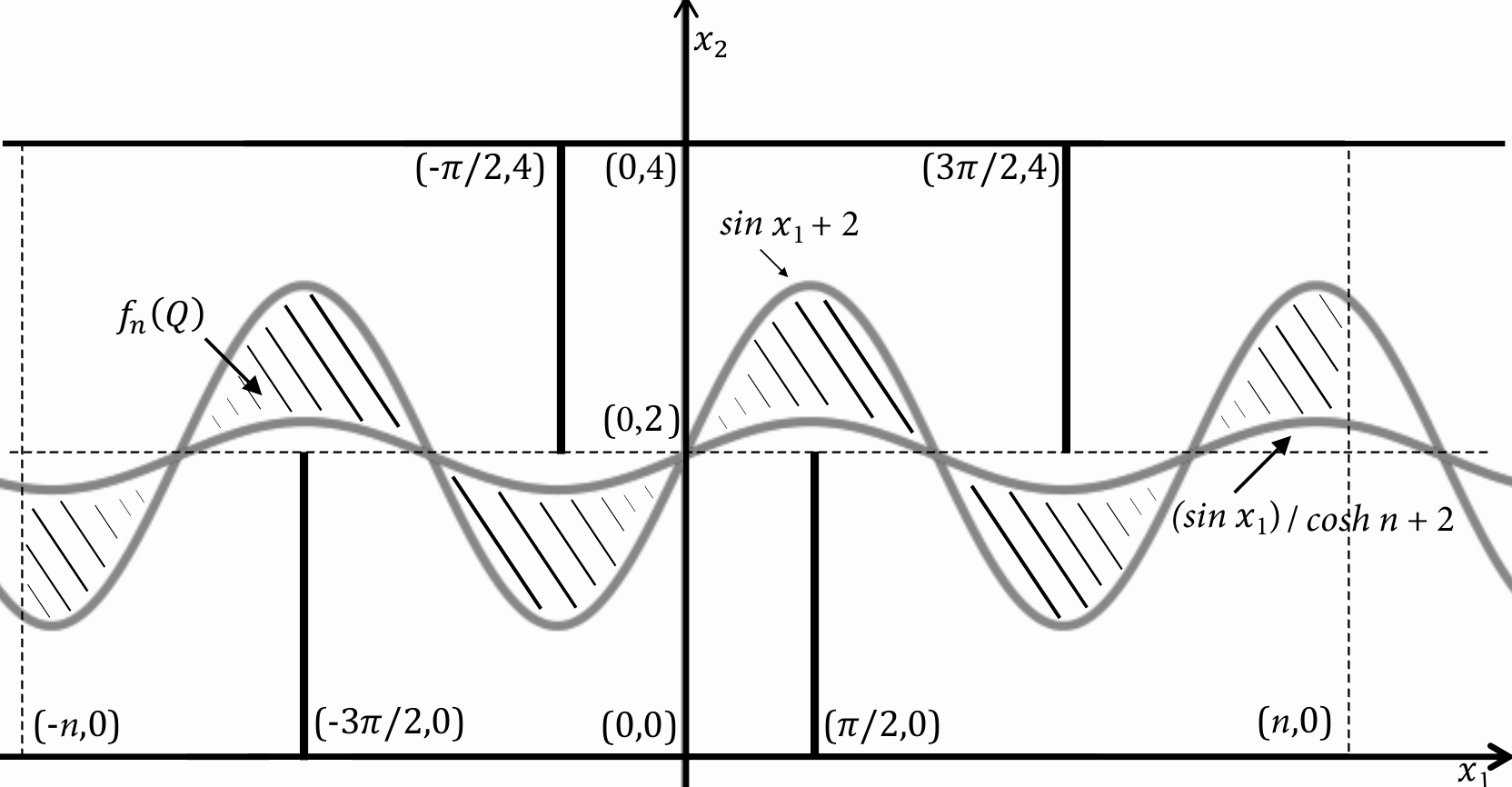}
\caption{\small An example with rough boundary.}  
\end{center}  
\end{figure} 

We will now prove that $T_D$ is not hyperbolic. Let $a:=(0,2)\in D$. We will construct a sequence of harmonic mappings $f_n:\Delta\to D$ such that $f_n(0)=a$ and $||df_n(0)||\ra\infty$ as $n\ra \infty$. Define
\begin{equation}
f_n:\Delta\to\RR^2,\quad z=x+iy\mapsto \left(nx, \frac{\sin(nx)\cosh(ny)}{\cosh n}+2\right),\quad n\in\NN.\label{mapf}
\end{equation}
Clearly, $f_n$ is harmonic, $f_n(0)=a$ and $||df_n(0)||\ra\infty$ as $n\ra \infty$.

It remains to see that $f_n$ maps $\Delta$ into $D$ for all $n\in\NN$. Consider the closed unit square in $\CC$
$$
Q:=\{z=x+iy\in\CC\mid -1\le x\le 1,\, -1\le y\le 1\}.
$$  
For every fixed $-1\le y_0\le 1$, the map $f_n$ takes the segment $[-1,1]\times\{y_0\}\subset Q$ into the curve
$$
\Gamma_{y_0}:=\left\{\left(x_1,\frac{\sin x_1\cosh(ny_0)}{\cosh n}+2\right)\biggm| -n\le x_1\le n\right\}.
$$
Therefore, the image
$$
f_n(Q)=\bigcup_{-1\le y_0\le 1}\Gamma_{y_0}
$$
is the closed set bounded by the graphs of $\sin x_1+2$ and $(\sin x_1)/\cosh n+2$ on the segment $-n\le x_1\le n$, which is shown as the shaded area in Fig.~1. Thus, $f_n(Q)$ lies in $D$ and so does $f_n(\Delta)\subset f_n(Q)$.

The above example can be modified by choosing
$$
D:=\{(x_1,x_2)\in\RR^2\mid 0<x_2<4\}\setminus (S_1\cup S_2\cup S_3\cup S_4),
$$
where each of 
$$
\begin{array}{l}
S_1\subset \{(x_1,x_2)\in\RR^2\mid x_1\le 0,\, 0\le x_2\le 2 \},\\
\vspace{-0.1cm}\\
S_3\subset \{(x_1,x_2)\in\RR^2\mid x_1\ge 0,\, 0\le x_2\le 2 \}
\end{array}
$$
is a closed region whose boundary contains a curve joining a pair of points on the line $\{x_2=0\}$ and passing through a point on the line $\{x_2=2\}$ and each of    
$$
\begin{array}{l}
S_2\subset \{(x_1,x_2)\in\RR^2\mid x_1\le 0,\, 2\le x_2\le 4 \},\\
\vspace{-0.1cm}\\
S_4\subset \{(x_1,x_2)\in\RR^2\mid x_1\ge 0,\, 2\le x_2\le 4 \}
\end{array}
$$
is a closed region whose boundary contains a curve joining a pair of points on the line $\{x_2=4\}$ and passing through a point on the line $\{x_2=2\}$, as shown in Fig. 2 below. Furthermore, we require that none of the regions $S_j$ intersects the graph of the function $\sin x_1+2$.

\begin{figure}[H]
\begin{center}  
\includegraphics[width=4.2in]{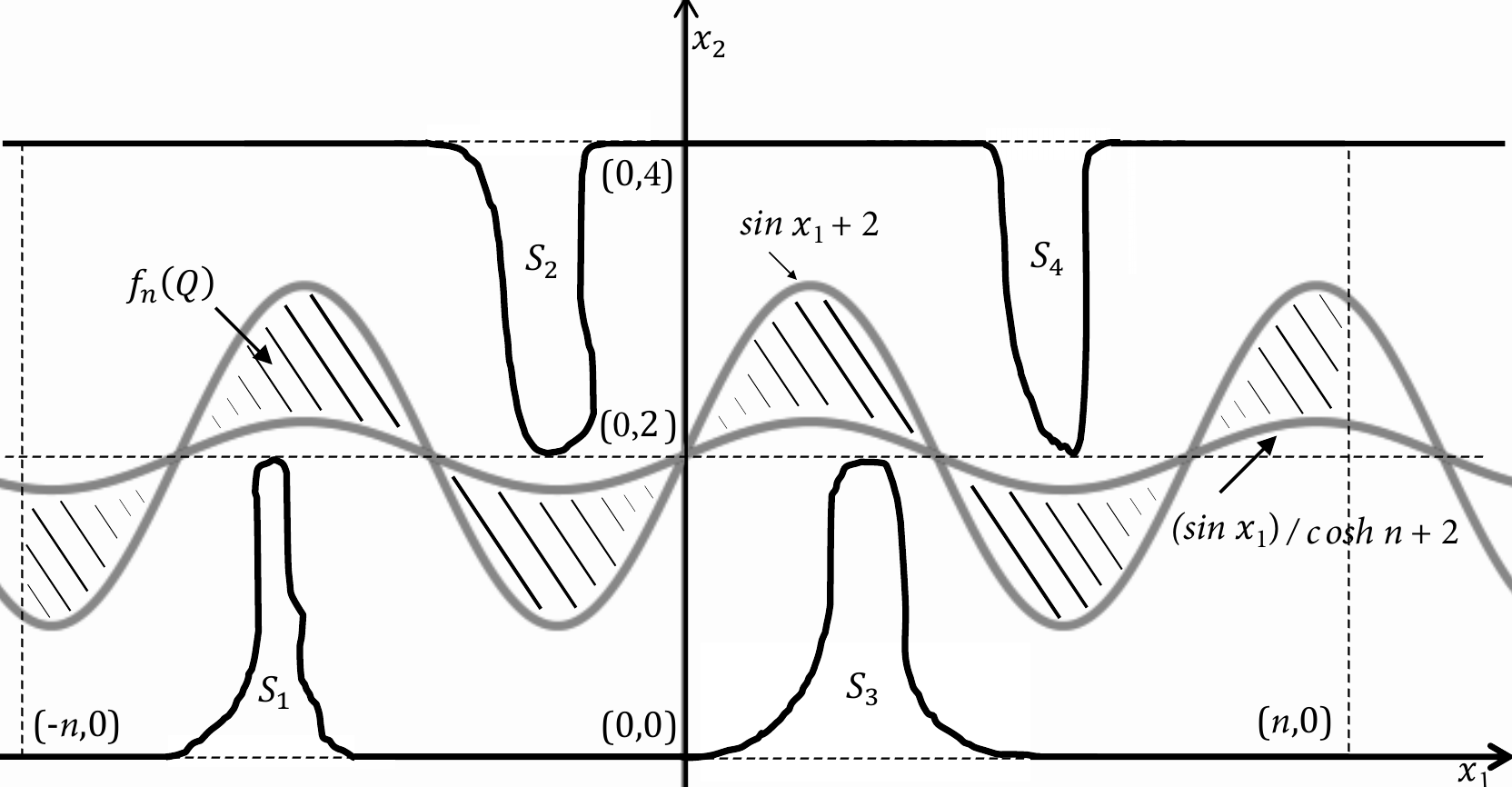}
\caption{\small An example with smooth boundary. }  
\end{center}  
\end{figure} 

\noindent Clearly, the regions $S_j$ can be chosen to ensure that $\partial D$ is smooth. 

As above, we now see that $D$ has Property (J-P$)_{\aff}$ and that $T_D$ is not hyperbolic. This yields the second statement of Theorem \ref{result}.\qed 

\begin{remark}\label{gamma} Observe that in the examples given above Property (J-P) fails for $a=(0,2)$ and
$$
\gamma_k(t)=\frac{\sin t}{k}+2,\quad t\in[-k,k].
$$
\end{remark}

\begin{remark}\label{holommap} It is easy to find harmonic conjugates to the components of the harmonic maps $f_n$ defined in (\ref{mapf}), which yields a sequence of holomorphic maps $g_n:\Delta\to T_D$ with $g_n(0)=a$ and $||dg_n(0)||\to\infty$ as $n\to\infty$:
$$
g_n(z):=\left(nz, \frac{\sin(nx)\cosh(ny)+i\cos(nx)\sinh(ny)}{\cosh n}+2\right)=\left(nz, \frac{\sin(nz)}{\cosh n}+2\right).
$$
\end{remark}

\end{document}